\documentclass[12pt]{article}
\usepackage[T1]{fontenc}
\usepackage{kpfonts}
\usepackage[dvipsnames]{xcolor}
\usepackage{pdfpages}
\usepackage{sgame}
\usepackage{accents}
\usepackage{tikz-cd}
\usepackage{float}
\usepackage[round]{natbib}
\usepackage{multirow}
\usepackage{dutchcal}
\usepackage{pdfpages}
\usepackage{bbm}
\usepackage{sgame}
\usepackage{mathtools}
\usepackage{amsthm}
\usepackage{enumitem}
\usepackage[margin=1in]{geometry}
\usepackage{caption}
\usepackage{subcaption}
\usepackage{epigraph}
\usepackage{listings}
\usetikzlibrary{calc}
\usetikzlibrary{shapes,arrows}
\usepackage{tcolorbox}

\newcommand{\ubar}[1]{\underaccent{\bar}{#1}}

\newtheorem{theorem}{Theorem}[section]
\newtheorem{proposition}[theorem]{Proposition}

\newtheorem{corollary}[theorem]{Corollary}

\theoremstyle{definition}

\newtheorem{remark}{Remark}

\usepackage{hyperref}
\hypersetup{
    colorlinks=true,
    linkcolor=ForestGreen,
    filecolor=Periwinkle,      
    urlcolor=Periwinkle,
    citecolor=Aquamarine,
}

\definecolor{backcolour}{rgb}{0.63, 0.79, 0.95}

\lstdefinestyle{mystyle}{
  backgroundcolor=\color{backcolour},
  basicstyle=\ttfamily\footnotesize,
  breakatwhitespace=false,         
  breaklines=true,                 
  captionpos=b,                    
  keepspaces=true,                 
  numbers=left,                    
  numbersep=5pt,                  
  showspaces=false,                
  showstringspaces=false,
  showtabs=false,                  
  tabsize=2
}

\providecommand{\keywords}[1]{\textbf{\textit{Keywords:}} #1}
\providecommand{\jel}[1]{\textbf{\textit{JEL Classifications:}} #1}

\begin{document}
\title{Dynamic Competitive Persuasion}
\author{Mark Whitmeyer\thanks{Arizona State University. Email: \href{mailto:mark.whitmeyer@gmail.com}{mark.whitmeyer@gmail.com}. I began this paper, then called ``A Game of Martingales," at UT Austin and continued to work on it while I was at the University of Bonn. The DFG (project ID 390685813) supported me during this period. I am grateful to Isaac M. Sonin for suggesting this project and to V. Bhaskar, Raphael Boleslavsky, Andy Kleiner, Joseph Whitmeyer, Thomas Wiseman, and Kun Zhang for their feedback. I also thank the editor, Tilman B\"orgers, and two anonymous referees for their advice.}}

\date{\today}

\maketitle

\begin{abstract}
Two long-lived senders play a dynamic game of competitive persuasion. Each period, each provides information to a single short-lived receiver. When the senders also set prices, we unearth a folk theorem: if they are sufficiently patient, virtually any vector of feasible and individually rational payoffs can be sustained in a subgame perfect equilibrium. Without price-setting, there is a unique subgame perfect equilibrium. In it, patient senders provide less information--maximally patient ones none.
\end{abstract}
\\
\jel{C72; C73; D15; D63}

\newpage

\section{Introduction}

How does competition over the long run shape information provision? Two long-lived senders interact over an infinite time horizon and compete each period through information design. Each period, each sender chooses a Blackwell experiment (or signal) about her independent binary type with the aim of attracting a sequence of short-lived receivers.

We analyze two cases. In the first, the senders compete through information alone, and each period's receiver selects the sender whose expected type is the highest. In the second case, the senders are sellers who also post prices each period. There, a receiver (consumer) purchases from a seller if and only if its expected type net of its price exceeds that of the other seller and the consumer's outside option (normalized to \(0\)).

The two setups (with and without prices) are spartan, yet include the following essential ingredients. There is competition, as only one sender is selected each period, to the detriment of the other sender. There is a tension between immediate and distant incentives, as the short-run objective of besting one's opponent in the current period's competition is tempered by the long-run ramifications of the current period's persuasion.\footnote{The state is persistent and the senders' memories long, so any information provided today will affect play forever.} There is also--when the senders also set prices--an extraction motive: there is no longer a fixed pie to split. Is the repeated interaction between senders to the receivers' benefit? How does the price instrument alter things? What is the role of the senders' patience?

The environment with prices captures competition between two duopolists who not only charge prices but provide information to consumers about their products. The setting without prices also corresponds to some stylized cases itself. Consider, for example, two universities competing to place a representative student from their graduating classes in a prestigious job. Each year, the schools dictate their grading policies to their professors--they commit to a signal about the hidden state, the quality of the student. Year-over-year, the pool of students that comprises each school's student body stays roughly the same, which generates the sort of gradual learning present in our model. A second environment that fits is a political one. Political opponents frequently run against each other multiple times throughout their careers and choose how much to reveal about their positions. Our stipulation that memories are long rings especially true here: there are quite a few infamous examples of past remarks (information revelation) coming back to wreak havoc on a candidate's campaign.

Information revelation (or risk-taking) comes at a cost: the history of experiments, and their realizations (and prices) is public, which leads to a (martingale) sequence of beliefs. When the senders do not post prices, we derive the unique subgame perfect equilibrium of the game. The game is formally similar to an all-pay auction; and, just as participants in such auctions typically choose continuous distributions over bids, the senders here choose continuous distributions over their ``bids,'' the posteriors that their experiments produce. Moreover, in contrast to the all-pay auction, where winning bids are more costly the higher they are, the senders here do not care about the realized posteriors, \textit{ceteris paribus}, so the senders choose distributions over posteriors that yield the other a \textit{linear} payoff. As in an all-pay auction, this linearity ensures that the other is indifferent over the distribution she chooses.

In a single period game of competitive persuasion, \cite{Au2} find that competition increases information provision: as the number of senders increases, each provides more and more information; in the limit, full information. Here, (absent prices) we discover that patience has the opposite effect: as senders become more patient, each provides less information. As they become maximally patient, the sender with the highest expected type provides no information. The other sender chooses a binary distribution with support on \(0\) and the mean of the high sender, then provides no information forever onward, herself. 

Information provision is a blunt tool: today's posterior is tomorrow's prior and any revelation constrains what can be done forever onward. Moreover, the game is constant-sum, so there is no scope for cooperation; thus, the two senders' equilibrium values are unique. As a result, the analysis of the dynamic game can be reduced to the study of a single period interaction in which each sender's payoff is today's payoff plus the expected (discounted) value of the game for that sender. As noted above, these payoffs must be linear in the realized posteriors and also yield the equilibrium values when evaluated at the vector of priors. The resulting pair of differential equations pins down the unique equilibrium distributions chosen by each sender.

Why do senders provide less information as they become more patient? The driving force behind this is the concavity of a sender's equilibrium value in her prior. These values (which are the senders' continuation values) become more important--and so payoffs become more concave--as the senders grow patient. A sender's equilibrium value is concave in her prior due to the fact that the richness of the strategy space--a sender can choose any Blackwell experiment--allows the disadvantaged sender (the one with the lower prior) to partially imitate the advantaged sender. This ensures diminishing returns to being advantaged: a disadvantaged sender's value is linear in her prior, and an advantaged sender's value is linear in the inverse of her prior, ergo strictly concave.


Allowing the senders to post prices grants them a great deal of flexibility and transforms the game into one that is no longer constant sum. They can adjust prices freely, which enables them to reward and punish each other as needed, essential to sustaining cooperation. The latitude granted by price setting engenders a folk theorem: for any vector of average payoffs that is in the set of feasible and strictly individually rational payoff vectors, there exists a subgame perfect equilibrium of the game that yields such payoffs provided the sellers are sufficiently patient. In contrast to the information-only game, in which long-run incentives dampen information provision, we obtain the result in the price-setting scenario by constructing equilibria in which the firms provide full information immediately before cooperating in the resulting repeated Bertrand game.

\subsection{Related Literature}

There are a number of papers that study the problem of competitive persuasion (without prices) in a single period setting. They are \cite{cotton}, who solve this problem for two players; \cite{Albrecht}, who sets the two-player competition in a political realm; \cite{Au2}, who extend the analysis to $n$ players; and \cite{Au1}, who allow for correlated types. \cite{spiegler} solves a mathematically identical problem of firms competing for boundedly rational consumers. 

Quite a few other works explore risk-taking contests in which agents with a fixed budget (or starting point) can ``gamble'' by choosing distributions on the positive portion of the real line. This problem is essentially equivalent to the competitive persuasion problem, with the minor difference that players choose distributions on $\left[0,\infty\right)$ rather than $\left[0,1\right]$. \cite{wagman2012choosing} and \cite{fang} focus on the martingale case, in which the expectations of the players' distributions must be their budgets. In \cite{seel2013continuous} the players choose when to stop Brownian motions with drift, and \cite{feng2015gambling} extend this to general diffusions.

The common theme throughout each competitive persuasion or gambling paper is the linearity of the equilibrium payoff function: players choose distributions that yield a payoff function for each other player that is linear in that player's realization, leaving that player, herself, willing to choose an equilibrium distribution. As we discover in here, the necessity of the linear payoff persists. However, this payoff now includes the continuation payoff, which becomes increasingly important as players become more patient.

There are also a few recent investigations of dynamic information disclosure by a single sender. \cite{audynamic} extends \cite{gent} to such a setting: a sender wishes an agent to take one irreversible action and discloses information over time. Both players are long-lived, and as the players become maximally patient, full disclosure occurs virtually immediately. This contrasts starkly with our finding that maximally patient players are completely reticent. Also related is \cite{guo}, who generalize \cite{audynamic} and allow for correlated sender and receiver types. 

A recent paper, \cite*{koes}, explores ``long information design" in which persuaders disclose information at multiple stages about their (respective) states before the receiver chooses an action. As in this paper, they find that less information is disclosed when there are multiple stages of disclosure than in the one stage (single period) case. In addition, the equilibrium of this paper's game coincides with the equilibrium in their product demonstration example when there is no deadline. 

In contrast to the abundance of pure persuasion papers, only a few study competition in which sellers post prices and choose information structures. \cite{gar}, \cite*{hkb}, and \cite*{boleslavsky2017demonstrations} all look at variations of such a problem;\footnote{The first two of these papers focus on symmetric settings and have the same timing as the stage-games of this paper: prices and information provision are determined simultaneously. In contrast, \cite{boleslavsky2017demonstrations} compare two different timings. In one, a firm chooses its information policy first (only one firm provides information in their duopoly model), which is observed before firms compete over prices. In the other, prices are chosen first, before the information policy is specified.} and \cite{whau} and \cite{dia} investigate the effect of search frictions on this type of competition. The infinite horizon of this paper alters this competition drastically--just as long-run interaction and patience facilitate cooperation (and surplus extraction) in the standard model of Bertrand competition, so do they here.

\section{The Game Without Prices}

Time is discrete, indexed by the natural numbers, $t = 0, 1, 2, \dots$. There are two long-lived senders and a sequence of short-lived receivers, who have identical preferences. Each sender has a persistent binary type taking values $\left\{0,1\right\}$ that is drawn prior to period \(0\) by nature and remains unchanged for the duration of the game. The senders' types are independent. \(x_{0}\) (\(y_{0}\)) denotes the prior probability at period \(0\) that sender \(1\) (\(2\)) is value \(1\).

Each period, without knowing her (or the other sender's) type, each sender chooses a public signal about her type. Following the realizations of the two signals, that period's receiver selects the sender with the highest expected type (and randomizes fairly between the two if indifferent). We normalize the reward to a sender from being chosen to \(1\), and to \(0\) from not being chosen. The two senders discount the future by a common discount factor $\delta \in \left[0,1\right)$. Each sender and receiver has perfect recall and observes the complete history of signals and signal realizations chosen by each sender.

Each sender has a compact metric space of signal realizations and each period $t$ independently chooses a signal. As is now ubiquitous in this literature, we immediately reframe this problem as one of senders choosing Bayes-plausible distributions over posterior values: each period, given priors $\left(x_t, y_t\right)$, senders \(1\) and \(2\) choose random variables $X_t \sim F_t$ and $Y_t \sim G_t$ supported on $\left[0,1\right]$ whose expected values are $x_t$ and $y_t$, respectively.\footnote{It is without loss of generality to restrict the senders to pure strategies, since any randomization over distributions is strategically equivalent to the resulting marginalized distribution.} Naturally, a sender's realized posterior in period $t$ becomes the prior in period $t+1$. $\mathcal{F}_\omega$ ($\omega = x, y$) denotes the set of Bayes-plausible distributions over posteriors.

\begin{remark}
    It is without loss of generality to identify signals with Bayes-plausible distributions over posteriors.
\end{remark}
\begin{proof}
    Any sequence of signals induces a sequence of distributions over posterior beliefs. Consequently, we need only to argue that we do not ``lose'' equilibria by working directly with the sequence of beliefs. But this is nearly immediate: the game at hand is constant-sum, and in the proof of Proposition \ref{spepayoffs}, we derive the two senders' equilibrium values, which depend only on the vector of prior means. This means that the dynamic game can be reduced to a single-period game, in which each sender's payoff is, moreover, posterior-separable--as today's payoff, and the discounted continuation payoff (which is the expected discounted value of the game) depend only on the realized posterior today. This posterior-separability implies the remark.\end{proof}

\subsection{Analysis of the Game without Prices}

We begin by characterizing each sender's equilibrium payoff. In the single period game of \cite{cotton}, the sender with the higher mean (WLOG sender \(1\)) obtains payoff $1 - \frac{y_0}{2x_0}$ and sender \(2\) gets $\frac{y_0}{2x_0}$. Here, as we discover, in any subgame perfect equilibrium of the infinite-horizon game, these are the two senders' average payoffs.
\begin{proposition}\label{spepayoffs}
In any subgame perfect equilibrium of the infinite horizon game, the two senders' average payoffs are their payoffs of the one-shot game.
\end{proposition}
The proof of this result may be found in the appendix but here is a sketch. It is clear that sender \(2\)'s equilibrium payoff must be at least \(\frac{1}{1-\delta}\left(\frac{y_0}{2x_0}\right)\): given any strategy by sender \(1\), sender \(2\) could guarantee herself this payoff by choosing a compound distribution in which she realizes \(0\) or \(x_{0}\), then, following the latter realization, mimics sender \(1\) thenceforth. However, it is \textit{a priori} unclear whether sender \(2\) can be limited to just this payoff--the proposition verifies that this is indeed the case: we construct a strategy for sender \(1\) that limits sender \(2\) to this value. As the game is constant-sum, this pins down the unique equilibrium payoffs for both senders.

Proposition \ref{spepayoffs} allows us to reduce the dynamic game to a single-period game. In any period, given any distribution \(H_2\) chosen by sender \(2\), sender \(1\) solves
\[\sup_{F \in \mathcal{F}_{x}}\left\{\int_{0}^{1}\Gamma_1\left(w\right)dF\left(w\right)\right\} \text{ ,}\]
where
\[\Gamma_1\left(w\right) \coloneqq H_2\left(w\right) + \frac{\delta}{1-\delta}\int_{0}^{w}\left(1-\frac{z}{2w}\right)dH_2\left(z\right) + \frac{\delta}{1-\delta}\int_{w}^{1}\frac{w}{2z}dH_2\left(z\right) \text{ .}\]
From here the rest of our analysis is easy. \cite{cotton} (as do, e.g., \cite{fang} and \cite{Au2}) identify that the equilibrium payoff function must have a linear form.\footnote{Refer to either Lemma 4.1 in \cite{cotton}, the supplementary appendix in \cite{fang} or Theorem 1 of \cite{Au2}.} There can be no jumps--corresponding to interior atoms--in a sender's payoff, since in her best response she prefers to place a mass just above the jump in her payoff; which, in turn, means that her opponent is not best responding. The support of the two senders' equilibrium distributions must be the same (or else one would prefer to spread out the measure on the subset uninhabited by her competitor), nor can there be gaps in the support except possibly at the top (or else senders would want to place mass points). Similar arguments preclude $\Gamma_1$ (or sender \(2\)'s payoff, $\Gamma_2$) from being strictly convex or concave.

Given the linear form, it remains to solve the differential equation(s) $\Gamma_i = \lambda_i + \eta_i w$ ($i = 1, 2$), accompanied by the conditions $\int_{0}^{1}wdF\left(w\right) = x$, $\int_{0}^{1}wdG\left(w\right) = y$, and that \(F\) and \(G\) are cdfs supported on (a subset of) the unit interval. This pins down the unique equilibrium distributions at every history. Defining $\mu \coloneqq \left(1 + \sqrt{1-\delta}\right)^{-1}$, the sender with the higher current mean (WLOG sender \(1\)) chooses a continuous, atomless, distribution on $\left[0, \frac{x}{\mu}\right]$ if \(x \leq \mu\) and a distribution with a continuous, atomless, portion on $\frac{1-x}{\left(1-\delta x\right)\mu}$ and a mass point on \(1\) if \(x \geq \mu\). sender \(2\) partially mimics sender \(1\): she chooses sender \(1\)'s distribution with probability $y/x$ and places a mass point on \(0\) with its complement (which is needed in order to satisfy the Bayes-plausibility condition). Summing things up, and leaving the details of the distributions and payoff functions to the appendix,
\begin{theorem}\label{mainthm}
There is a unique subgame-perfect equilibrium of the game. It is Markovian, in which the state variable $(x,y)$ is the current vector of priors. Each sender chooses a distribution that generates a payoff function for her opponent that is linear in the opponent's posterior.
\end{theorem}
With Proposition \ref{spepayoffs} in hand, Theorem \ref{mainthm} is easy, and the economics behind the result is straightforward. In short, it is a standard ``matching pennies'' argument: each sender's continuous distribution must give the other a linear payoff so that she, herself, is willing to choose a continuous distribution. The necessity of this linear form then makes it clear why the proposition holds: the senders are indifferent over any distribution supported on $\left[0,\frac{x}{\mu}\right]$, two of which are the degenerate distribution supported on $\left\{x\right\}$ (for sender \(1\)) and the distribution for sender \(2\) that (partially) mimics this, i.e., the distribution supported on $\left\{0,x\right\}$ with mean $y$.

As senders become more patient, the relative importance of besting one's opponent decreases. Indeed, consider first the immediate payoff (whose relative magnitude wanes with patience): for a fixed realization of a sender's opponent, $z$, it is a step function that jumps from $0$ to $1$ as soon as the sender's posterior exceeds $z$. On the other hand, the long-run payoff exhibits no such jump: for a fixed realization $z$ it is linearly increasing in the sender's posterior on $\left[0,z\right)$, then strictly concave on $\left[z, 1\right]$; and is, therefore, concave on the entire unit interval. Moreover, this payoff function is not only continuous on \(\left(0,1\right]\) but differentiable. These observations suggest the next result: today's payoff encourages information provision, as each sender wants to provide a marginally higher posterior than her opponent; whereas the continuation payoff discourages such disclosure, due to its inherent concavity. Consequently, as the latter becomes more important, senders should provide less information.

The previous paragraph's intuition was mostly mechanical. It argues that the decrease in information provision as senders become more patient is driven by the increased importance of the continuation value, which is concave. But why is it concave? The fact that at equilibrium, the disadvantaged sender does a mimicking strategy--whereby its distribution corresponds to an initial binary distribution with support on \(0\) and the advantaged sender's mean, the latter of which is followed by an exact duplicate of the advantaged sender's distribution--means that the disadvantaged sender's equilibrium payoff is linear in its mean. Then, the constant-sum aspect of the game means that the advantaged sender's payoff is linear in the inverse of her mean; and is, therefore, concave. In short, the ability of the disadvantaged sender to imitate means that being disadvantaged is not so bad and that there are diminishing returns to being advantaged.

\begin{proposition}\label{risky}
As senders become more patient they provide less information. That is, if $\hat{\delta} \geq \delta$ then the equilibrium distributions for discount factor $\delta$ are mean-preserving spreads of those distributions for discount factor $\hat{\delta}$.
\end{proposition}

As the senders become maximally patient, the advantaged sender becomes completely uninformative and the disadvantaged sender's distribution converges to a binary distribution where one of the realizations leaves the receiver just willing to select her and the other realization reveals that her type is \(0\). To wit, 
\begin{corollary}\label{koes}
As $\delta \to 1$, the equilibrium distributions converge to the following: sender \(1\) chooses the degenerate distribution with support $\left\{x\right\}$ and sender \(2\) chooses the binary distribution with support $\left\{0,x\right\}$.
\end{corollary}
In the limit, today's payoff is negligible, and all that matters is the continuation payoff, which is concave for any fixed realization of a sender's opponent. Consequently, the advantaged sender prefers not to release any information, whereas the disadvantaged sender's payoff is linear on $\left[0,x\right]$, leaving her willing to choose the specified binary distribution.

As noted in the introduction, the equilibrium described in Corollary \ref{koes} coincides with the equilibrium in the product demonstration application (with no deadline) of \cite{koes}. This equivalence is implied by their Lemmas 9 and 10 (which follow from results in \cite{laraki}), which establishes that their demonstration game without a deadline can be approximated by a discounted (splitting) game as senders become maximally patient.

An earlier version of this paper focuses on a finite-horizon (\(T\)-period) model.\footnote{Please refer to v4 of this paper on ArXiv.} In the unique subgame-perfect equilibrium of that game, obtained via backward induction, the senders' strategies and equilibrium values are the natural analogs of their equilibrium strategies and values here. Notably, as the horizon grows long (\(T \uparrow \infty\)), the equilibrium converges to the one identified in Theorem \ref{mainthm}.

\section{The Game With Prices}

We carry over several of the assumptions from the game without prices. Namely, time is discrete and the game is played between two long-lived senders, who wish to entice a sequence of short-lived receivers. Now, the two senders are sellers and the receivers consumers. Each period each seller publicly chooses a Blackwell experiment \textit{and} posts a price. Given a price $p$ and realized value from a seller's good $z$, a consumer's utility from purchasing from that seller is $z - p$. For simplicity, we stipulate that the consumer's outside option is \(0\).

No longer do we limit the analysis to binary types: initially, at time $t=0$, the values of seller \(1\)'s and seller \(2\)'s goods to the consumer are independent random variables $X_0 \sim F_{0}$ and $Y_0 \sim G_{0}$, respectively, that are each supported on subsets of the unit interval. In principle, this infinite-horizon game of price and information competition is a very difficult problem: each period's information provision decision by a firm determines tomorrow's prior, and because of this the firms' decisions do not reduce to choosing mean-preserving contractions of the prior, as they would if this were a single-period problem. However, as it turns out we can bypass these difficulties entirely.

Before doing so, let us reflect on a few things. If this were a market with just one seller, the optimal pricing and information policy for the (now) monopolist would be to never provide any information and charge a price equal to the good's expected quality. Trade is efficient--since the consumer always buys--and the seller extracts all of the surplus. Furthermore, if this monopolist were extremely patient, it would be almost optimal for her to immediately provide full information, then set a price equal to the consumer's value--as the discount factor approaches \(1\) the seller would approximate her optimal payoff. 

On the other hand, with two sellers, it is not efficient for the sellers to provide no information. With two sellers, efficiency requires not only that the consumer always purchase from a seller (since her outside option is \(0\)) but also that the consumer purchase from the seller with the higher value good. Accordingly, it is efficient for both sellers to provide full information; and, therefore, if the sellers are extremely patient, it is approximately optimal for the sellers to provide full information immediately and for the higher value seller to sell to the consumer at a price equal to the consumer's value.

The (near-)optimality of immediate information revelation suggests an approach in the repeated game between the sellers. We should have the sellers immediately (in period \(0\)) provide full-information--maximizing the available surplus--then cooperate over the remainder of the interaction in splitting the surplus. As it turns out, this also helps us in constructing the equilibrium: from $t=1$ on, (on path) there is no longer uncertainty and commencing in any such subgame, the sellers are now in a standard repeated game. Moreover, it is easier to tackle unilateral deviations in the initial period, since the other (non-deviating) seller will have a much simpler strategy space from $t=1$ on. Thenceforth, her prior is degenerate and therefore she can only affect payoffs through the price instrument.

Given distributions over valuations \(F\) and \(G\), the maximal surplus is
\[S_{F,G} \coloneqq \mathbb{E}\left[\max\left\{X,Y\right\}\right]\text{ .}\]
Understanding $V_{i, min} = V_{i, min}\left(F,G\right)$ and $V_{i, max} = V_{i, max}\left(F,G\right)$ for $i = 1,2$, we define \[V_{1, min} \coloneqq \int_{0}^{1}\int_{0}^{x}\left(x-y\right)dG\left(y\right)dF\left(x\right), \quad V_{2, min} \coloneqq \int_{0}^{1}\int_{0}^{y}\left(y-x\right)dF\left(x\right)dG\left(y\right) \text{ ,}\]
\[V_{1, max} \coloneqq \int_{0}^{1}xdF\left(x\right), \quad V_{2, max} \coloneqq \int_{0}^{1}ydG\left(y\right) \text{ ,}\]
and
\[V\left(F, G\right) \coloneqq \left\{\left(v_1, v_2\right) \in \left[0,1\right]^2 \colon V_{1, min} < v_1 < V_{1, max} \ \text{\&} \ V_{2, min} < v_2 < V_{2, max} \ \text{\&} \ v_1 + v_2 < S_{F, G} \right\} \text{ .}\]
$V_{i,min}$ is the average payoff a seller $i$ can (approximately) guarantee herself when she is very patient. $V_{i,max}$ is seller $i$'s monopoly payoff, which is obviously the maximal payoff seller $i$ can obtain. $V\left(F,G\right)$ is our analog of the feasible and individually rational set.
\begin{theorem}\label{pricetheorem}
For any $\left(v_1,v_2\right) \in V\left(F_0,G_0\right)$ there exists a discount factor $\ubar{\delta} < 1$ such that for any $\delta \geq \ubar{\delta}$ there is a Subgame Perfect Equilibrium of the game in which the average payoff to seller $i$ is $v_i$.\end{theorem}
The detailed construction and proof of this result may be found in the appendix. We assemble the equilibria by having both sellers immediately provide full information, then cooperate in one of the ensuing subgames.

Several more comments regarding Theorem \ref{pricetheorem} are in order. The theorem is an ``anything goes'' result for surplus splits, provided senders are sufficiently patient. It is not a complete characterization of equilibrium behavior. Indeed, it is easy to see that behavior other than the immediate full revelation we specified can be sustained at equilibrium. For example, provided senders are sufficiently patient, there are a variety of equilibria in which they provide less than full information in the first \(k\) periods then provide full information in the \(\left(k+1\right)\)th. In contrast, Theorem \ref{mainthm} is an exact characterization of the unique equilibrium strategies for arbitrary \(\delta\). 

The reason why we can obtain such a sharp characterization when there are no prices but instead must settle for a folk theorem when there are prices is because the former game is constant-sum and the latter game is not; and because the senders in the latter game have aspects of their action sets that are not persistent. Obviously, the constant-sum aspect of the game without prices pins down the unique values for the two senders, negating any multiplicity. With prices, the set of feasible and individually rational payoffs is no longer a singleton. Moreover, today's prices do not affect future prices, so senders have the requisite flexibility to cooperate.

\subsection{One-Shot Interaction}

Furthermore, \textit{both} the infinite horizon (and patience) and the price instrument are needed to generate the equilibrium multiplicity described in Theorem \ref{pricetheorem}. \cite{whau} study a single-period, binary-value, directed search model in which sellers compete by providing information and setting prices. Proposition 3.7 of that paper reveals that there is a unique symmetric equilibrium in which sellers provide full information and randomize over prices. The one-shot version of this paper's game is a special case of \citeauthor{whau}'s setting, when the search cost is \(0\). \cite{hkb} also identify a unique equilibrium in a game with price-setting and information provision. There, the distribution over posteriors has an alternating concave-linear form, where the nonlinearities arise due to constraints by the prior, and the prices are deterministic. 

\section{Discussion}

Absent prices, the senders in this paper face a challenging problem. Each period, they are tasked with designing an experiment to entice a receiver, yet it and its realizations are fully observable and so affect the future beliefs of each player. These dynamic concerns engender caution. As senders become more patient, they choose distributions that are more and more \textit{uninformative}.

Even though the game is infinite horizon, there is no scope for cooperation: any distribution that yields a sender a payoff that is not linear in her posterior gives her the opportunity to take (permanent) advantage of her opponent, and the constant-sum nature of the game means there can be no give and take between the senders. However, once price-setting is introduced, the permanence of information provision impedes cooperation no more. In fact, the sellers can extract maximal surplus by providing information immediately, then cooperating to eliminate the consumers' rents.

Our construction of equilibria when prices are present seems realistic. That is, we achieve cooperation by having firms reveal information immediately, then proceed with a price-setting game. Anecdotally, this seems to accurately reflect the typical behavior in a product's life cycle--there is an initial flurry of advertising for a new product, then little thereafter. Our result in the persuasion-only game also seems to correspond nicely to outcomes in the applications: for instance, the caginess of politicians is well-known, and schools, especially those whose input quality is highly esteemed, grade coarsely and/or provide little information about their students' performance.\footnote{A stark example of this is prohibition of grade disclosure by many top MBA programs--Chicago (Booth), Columbia, UC Berkeley (Haas), Stanford GSB, and UPenn (Wharton).}

\bibliography{sample.bib}

\appendix

\section{Omitted Proofs}

\subsection{Proposition \ref{spepayoffs} Proof}
\begin{proof}
The proof of this result is constructive. Recall the definition \[\mu = \mu\left(\delta\right) \coloneqq \frac{1}{1 + \sqrt{1-\delta}} \text{ .}\]
We also introduce the cdfs (with $l$ and $h$ representing ``low'' and ``high,'' respectively and $\xi \coloneqq \max\left\{x,y\right\}$).
\[\tag{$A1$}\label{eqa1}G^{\xi}_{l}\left(w\right) = \left(\frac{\mu w}{\xi}\right)^{\frac{1}{\sqrt{1-\delta}}} \quad \text{on} \quad \left[0,\frac{\xi}{\mu}\right] \text{ ,}\]
and
\[\tag{$A2$}\label{eqa2}G^{\xi}_{h}\left(w\right) = \begin{cases}
\left(1-a\right)\left(\frac{\left(1-\xi\delta\right)\mu w}{\left(1-\xi\right)}\right)^{\frac{1}{\sqrt{1-\delta}}} \quad &\text{if} \quad 0 \leq w < \frac{\left(1-\xi\right)}{\left(1-\xi\delta\right)\mu}\\
1-a, \quad &\text{if} \quad \frac{\left(1-\xi\right)}{\left(1-\xi\delta\right)\mu} \leq w < 1\\
1, \quad &\text{if} \quad  w \geq 1
\end{cases}\text{ ,}\]
where 
\[a = a\left(\xi\right) \coloneqq 2 - \frac{1}{\xi} - \left(\frac{\delta}{1-\delta}\right)\frac{\left(1-\xi\right)^{2}}{\xi} \text{ .}\]
Let sender \(2\) play the following Markovian strategy. Each period she chooses distribution \(G\) (though henceforth we suppress the subscript) defined as follows. When \(x \leq y\), \(G= G_l^y\) if \(y \leq \mu\) and \(G = G_h^y\) if \(y > \mu\). When \(y < x\), \(G = \frac{x}{y} G_l^x + 1 - \frac{x}{y}\) if \(x \leq \mu\) and \(G = \frac{x}{y} G_h^x + 1 - \frac{x}{y}\) if \(x > \mu\).

Given this, sender \(1\) faces a Markov decision problem. Her payoff is (recursively):
\[V\left(x, y\right) = \sup_{F \in \mathcal{F}_{x}}\left\{\int_{0}^{1}G\left(w\right)dF(w) + \delta\int_{0}^{1}\int_{0}^{1}V\left(w,z\right)g(z)dzdF(w)\right\} \text{ .}\]
Let us guess that the following strategy is an optimal strategy: if $x \geq y$, sender \(1\) chooses the degenerate distribution with support on $\left\{x\right\}$; and if $x < y$, she chooses a binary distribution with support $\left\{0,y\right\}$. A direct substitution and some algebra yields sender \(1\) a continuation payoff of $\frac{1}{1-\delta}\left(1-\frac{y}{2x}\right)$ for realizations $x \geq y$ and $\frac{1}{1-\delta}\frac{x}{2y}$ for realizations $y > x$. Appealing to the one-shot deviation principle, sender \(1\) solves
\[\sup_{F \in \mathcal{F}_{x}}\left\{\int_{0}^{1}\Phi\left(w\right)dF\left(x\right)\right\} \text{ ,}\]
where
\[\Phi\left(w\right) \coloneqq G\left(w\right) + \frac{\delta}{1-\delta}\int_{0}^{w}\left(1-\frac{z}{2w}\right)g(z)dz + \frac{\delta}{1-\delta}\int_{w}^{1}\frac{w}{2z}g(z)dz \text{ .}\]
Substituting in for \(G\), we have when $\mu \geq x \geq y$,
\[\Phi\left(w\right) = \begin{cases}
\frac{1}{1-\delta}\left[1+y\frac{w-2x}{2x^2}\right], \quad &\text{if} \quad 0 \leq w < \frac{x}{\mu}\\
\frac{1}{1-\delta}\left[1-\frac{y}{2x}\left(1-\sqrt{1-\delta}\right)\right], \quad &\text{if} \quad \frac{x}{\mu} \leq w \leq 1
\end{cases} \text{ ,}\]
when $\mu \geq y \geq x$,
\[\Phi\left(w\right) = \begin{cases}
\frac{1}{1-\delta}\left[\frac{w}{2y}\right], \quad &\text{if} \quad 0 \leq w < \frac{y}{\mu}\\
\frac{1}{1-\delta}\left[\frac{1}{2\mu}\right], \quad &\text{if} \quad \frac{y}{\mu} \leq w \leq 1
\end{cases} \text{ ,}\]
when $x \geq \max\left\{\mu,y\right\}$,
\[\Phi\left(w\right) = \begin{cases}
\frac{1}{1-\delta}\left[1+y\frac{w-2x}{2x^2}\right], \quad &\text{if} \quad 0 \leq w < \frac{\left(1-x\right)}{\left(1-x\delta\right)\mu}\\
\frac{1}{1-\delta}\left[1 + \frac{y\left(\frac{\left(1+\sqrt{1-\delta}\right)\left(1-x\right)}{1-\delta x}-2x\right)}{2x^2}\right], \quad &\text{if} \quad \frac{\left(1-x\right)}{\left(1-x\delta\right)\mu} \leq w < 1\\
\frac{1}{1-\delta}\left[\frac{y}{x}\left(1-\frac{a\left(x\right)}{2}-\delta \frac{x - a\left(x\right)}{2}\right)+1-\frac{y}{x}\right], \quad &\text{if} \quad 1 \leq w
\end{cases} \text{ ,}\]
and when $y \geq \max\left\{\mu, x\right\}$,
\[\Phi\left(w\right) = \begin{cases}
\frac{1}{1-\delta}\left[\frac{w}{2y}\right], \quad &\text{if} \quad 0 \leq w < \frac{\left(1-y\right)}{\left(1-y\delta\right)\mu}\\
\frac{1}{1-\delta}\left[\frac{\left(1-y\right)}{2y\left(1-y\delta\right)\mu}\right], \quad &\text{if} \quad \frac{\left(1-y\right)}{\left(1-y\delta\right)\mu} \leq w < 1\\
\frac{1}{1-\delta}\left[1-\frac{a\left(y\right)}{2}-\delta \frac{y - a\left(y\right)}{2}\right], \quad &\text{if} \quad 1 \leq w
\end{cases} \text{ ,}\]

Evidently, our specified strategy is indeed optimal. Thus, sender \(2\) can always guarantee herself the payoffs $\frac{1}{1-\delta}\left(1-\frac{x_{0}}{2y_{0}}\right)$ when $y_0 \geq x_0$ and $\frac{1}{1-\delta}\left(\frac{y_{0}}{2x_{0}}\right)$ when $y_0 \leq x_0$. Since sender \(1\) can do likewise, we conclude that the unique equilibrium payoff vector is as stated. \end{proof}

\subsection{Theorem \ref{mainthm} Proof}
\begin{proof}
The necessity of the linear payoff at equilibrium follows from \cite{cotton}. Sufficiency is immediate by concavification. Finally, it is straightforward to verify that the distribution \(G\) (Expressions \ref{eqa1} and \ref{eqa2}) specified in Proposition \ref{spepayoffs}'s proof is the unique solution to the differential equation that arises when specifying the linear payoff and imposing the various conditions (Bayes-plausibility and support). 
\end{proof}

\subsection{Proposition \ref{risky} Proof}
\begin{proof}
Fix means $y \leq x$. Let \(G\) be the equilibrium distribution chosen by seller \(1\) for discount factor $\delta$ and \(\hat{G}\) be its analog for discount factor $\hat{\delta} > \delta$. Let $x \leq \mu\left(\delta\right)$. It suffices to show that \(G\) crosses \(\hat{G}\) on \(\left(0,1\right)\) once from above on $\left(0,\frac{x}{\hat{\mu}}\right)$. Obviously $G\left(0\right) = \hat{G}\left(0\right)$. Moreover, $\frac{x}{\mu} > \frac{x}{\hat{\mu}}$ and so \(G\) and \(\hat{G}\) intersect at least once on $\left(0,\frac{x}{\hat{\mu}}\right)$. It remains to show that this intersection point is unique. Replacing $\sqrt{1-\delta}$ and $\sqrt{1-\hat{\delta}}$ with $\beta$ and $\hat{\beta}$, respectively, on $\left[0,x\left(1+\hat{\beta}\right)\right]$, we define
\[T\left(w\right) \coloneqq G\left(w\right) - \hat{G}\left(w\right) = \left(\frac{w}{x\left(1+\beta\right)}\right)^{\frac{1}{\beta}} - \left(\frac{w}{x\left(1+\hat{\beta}\right)}\right)^{\frac{1}{\hat{\beta}}} \text{ .}\]
Directly, \[T'\left(w\right) = \frac{1}{w}\left[\frac{1}{\beta}\left(\frac{w}{x\left(1+\beta\right)}\right)^{\frac{1}{\beta}} - \frac{1}{\hat{\beta}}\left(\frac{w}{x\left(1+\hat{\beta}\right)}\right)^{\frac{1}{\hat{\beta}}}\right]\text{ .}\]
Pick some $a \in \left(0,x\left(1+\hat{\beta}\right)\right)$ at which point $T\left(a\right) = 0$ (we argued above that such a point exists). Since $\beta > \hat{\beta}$, we have 
\[0 = \frac{T\left(a\right)}{\hat{\beta}a} = \frac{1}{\hat{\beta}a}\left[\left(\frac{a}{x\left(1+\beta\right)}\right)^{\frac{1}{\beta}} - \left(\frac{a}{x\left(1+\hat{\beta}\right)}\right)^{\frac{1}{\hat{\beta}}}\right] > \frac{1}{a}\left[\frac{1}{\beta}\left(\frac{a}{x\left(1+\beta\right)}\right)^{\frac{1}{\beta}} - \frac{1}{\hat{\beta}}\left(\frac{a}{x\left(1+\hat{\beta}\right)}\right)^{\frac{1}{\hat{\beta}}}\right] = T'\left(a\right) \text{ .}\]
Accordingly, whenever $T$ crosses the axis on $\left(0,x\left(1+\hat{\beta}\right)\right)$, it must do so from above, which can therefore happen only once. The same steps can be followed to show that \(G\) crosses \(\hat{G}\) once from above when $x \geq \mu\left(\delta\right)$.
\end{proof}

\subsection{Theorem \ref{pricetheorem} Proof}
\begin{proof}
At any $t$, the state vector is the pair of then-current prior distributions $\left(F_t, G_t\right)$. The following remark is useful. Suppose at some point each seller's quality has been fully revealed. Denote this period $t^{*}$--to clarify, this is the first period $t$ such that the prior held by a consumer for each seller is degenerate.\footnote{Shortly we will stipulate that $t^{*} = 1$ on path, i.e., it is only initially ($t = 0$) that a consumer's belief about a seller may be nondegenerate.} In that case, the state vector from period $t^{*}$ on is $\left(\delta_a, \delta_b\right)$, for some $\left(a, b\right) \in \left[0,1\right]^2$, where WLOG $a \geq b$. For any such vector, the maximal average total payoff for the two sellers is $a$. Then, in the parlance of \cite{folkthm} (henceforth \hypertarget{flk}{FM86}), the feasible and strictly individually rational set is
\[R\left(a,b\right) \coloneqq \left\{x, y \in \left[0,1\right]^2 \colon a-b < x \leq a \ \text{\&} \ 0 < y \leq a - x\right\} \text{ .}\]

By Theorem 1 in \hyperlink{flk}{FM86},
\begin{remark}\label{usefulremark}
In the subgame commencing at period $t^{*}$ in which the value realizations for sellers \(1\) and \(2\) are $a$ and $b$ ($\leq a$), respectively, for any $\left(x,y\right) \in R\left(a,b\right)$ there exists a discount factor $\ubar{\delta}$ such that for any $\delta \geq \ubar{\delta}$ there is a subgame perfect equilibrium in which the vector of average payoffs for the two firms is $\left(x, y\right)$.
\end{remark}
Given this, we construct the equilibrium by having the firms provide full information in the first period ($t=0$) and charge $p = 0$, then play the appropriate subgame perfect equilibrium in the (from period \(1\) on) \textit{repeated}--that is, no longer dynamic--game. This construction evidently yields the set $V\left(F_0,G_0\right)$. It remains to tackle deviations and specify strategies at subgames unreachable via unilateral deviations.

We can divide the remainder into the following three cases.
\begin{itemize}[noitemsep,topsep=0pt]
    \item \textcolor{OrangeRed}{Case 1:} histories in which the state vector at some $t \geq 1$ is $\left(\delta_{a},\delta_{b}\right)$ for some $\left(a, b\right) \in \left[0,1\right]^2$, and at least one sender deviated at some previous point.
    \item \textcolor{MidnightBlue}{Case 2:} histories in which the state vector at some $t \geq 1$ is $\left(F, G\right)$, where at least one of \(F\) or \(G\) is not degenerate, and in which both senders deviated previously.
    \item \textcolor{Violet}{Case 3:} histories in which the state vector at some $t \geq 1$ is $\left(F, \delta_b\right)$ ($b \in \left[0,1\right]$), where \(F\) is not degenerate, and in which only sender \(1\) deviated at some previous point.
\end{itemize}
Now let us specify strategies in each case. We impose that at least one of $F_0$ or $G_0$ is not degenerate (or else \hyperlink{flk}{FM86}'s folk theorem yields the result immediately), in which case $t^{*} \geq 1$. Define $p_i\left(t\right)$ to be \(1\) plus the number of times seller $i$ charged a non-zero price before period $t$.
\begin{tcolorbox}[colframe=OrangeRed,colback=white]\textbf{Case 1:}  WLOG $a \geq b$. We impose that if both senders deviated at some preceding $t' < t^*$ they play a SPE from $t^*$ on in which the vector of average payoffs is $\left(a-b +\frac{\eta_1^{p_{1}\left(t^{*}\right)}}{t^{*}}, \frac{\eta_{2}^{p_{2}\left(t^{*}\right)}}{t^{*}}\right)$ for some small $1 \gg \eta_1, \eta_2 > 0$. If only one sender deviated at some previous point, we stipulate that the vector of average payoffs from $t^*$ on is $\left(a-b +\frac{\varepsilon^{p_{1}\left(t^{*}\right)}}{t^{*}}, b -\frac{\varepsilon^{p_{1}\left(t^{*}\right)}}{t^{*}}\right)$, if it was sender \(1\) who deviated, and $\left(a-\frac{\varepsilon^{p_{2}\left(t^{*}\right)}}{t^{*}}, \frac{\varepsilon^{p_{2}\left(t^{*}\right)}}{t^{*}}\right)$, if it was sender \(2\) who deviated, for some small $1 \gg \varepsilon > 0$. \hyperlink{flk}{FM86} guarantees that such SPE exist (starting from period $t^*$) provided $\delta$ is sufficiently large.
\end{tcolorbox}
\begin{tcolorbox}[colframe=MidnightBlue,colback=white]\textbf{Case 2:} We specify that both senders charge price $p = 0$ and provide full information. Consequently, the state in $t+1$ is $\left(\delta_a,\delta_b\right)$ (for some $a, b \in \left[0,1\right]$), which is case 1, and so the vector of average payoffs from then on is either $\left(a-b +\frac{\eta_1^{p_1\left(t\right)}}{t+1}, \frac{\eta_{2}^{p_2\left(t\right)}}{t+1}\right)$ or $\left(\frac{\eta_1^{p_1\left(t\right)}}{t+1}, b-a+ \frac{\eta_{2}^{p_2\left(t\right)}}{t+1}\right)$, depending on whether $a$ or $b$ is higher, and where $1 \gg \eta_1, \eta_2 > 0$ are small. On path, $p_i\left(t\right) = p_i\left(t+1\right)$ for $i = 1, 2$.

The two senders' average payoffs are
\end{tcolorbox}
\begin{tcolorbox}[colframe=MidnightBlue,colback=white]
\[u_1 \coloneqq \delta \left\{\frac{\eta_1^{p_1\left(t\right)}}{t+1} + \int_{0}^{1}\int_{b}^{1}\left(a-b\right) dF_t\left(a\right)dG_t\left(b\right)\right\} \text{ ,}\]
and
\[u_2 \coloneqq \delta \left\{\frac{\eta_2^{p_2\left(t\right)}}{t+1} + \int_{0}^{1}\int_{0}^{b}\left(b-a\right) dF_t\left(a\right)dG_t\left(b\right)\right\} \text{ .}\]
Appealing to the one-shot deviation principle, suppose sender \(1\) deviates by charging $p' \in \left(0, 1\right]$ and choosing some signal.\footnote{It is clear that sender \(1\) cannot profit by a deviation in which she continues to charge price \(0\). This merely delays the cooperative (repeated-game) portion of the dynamic game.} sender \(1\)'s payoff is bounded above by
\[\left(1-\delta\right)\int_{0}^{1}\int_{b}^{1}\left(a-b\right) dF_t\left(a\right)dG_t\left(b\right) + \delta \left\{\frac{\eta_1^{p_1\left(t\right)+1}}{t+1} + \int_{0}^{1}\int_{b}^{1}\left(a-b\right) dF_t\left(a\right)dG_t\left(b\right)\right\} \text{ ,}\]
which is less than $u_1$ for all sufficiently large $\delta$. The same logic implies that sender \(2\) cannot deviate profitably to some $p'' > 0$. Any immediate gain will vanish as $\delta$ approaches \(1\) and her payoff in the cooperative portion will be discretely lower.
\end{tcolorbox}
\begin{tcolorbox}[colframe=Violet,colback=white]\textbf{Case 3:} This is analogous to the previous case: we specify that both senders charge $p = 0$ and that seller \(1\) provides full information. The remainder proceeds in a kindred manner \textit{mutatis mutandis}.
\end{tcolorbox}

Finally, we need to check incentives in period \(0\). WLOG we do so for seller \(1\). Recall that seller \(2\) is providing full information. On path, seller \(1\) obtains some $v_1 \in \left(V_{1,min}, V_{1,max}\right)$. If it deviates at $t = 0$, for any $\rho > 0$ there exists a $\hat{\delta} < 1$ such that for all $\delta \geq \hat{\delta}$, seller \(1\)'s average payoff is bounded above by \[\int_{0}^{1}\int_{y}^{1}\left(x-y\right) dF_0\left(x\right)dG_0 + \rho \text{ .}\]
Observing $V_{1,min}\left(F_0,G_0\right) = \int_{0}^{1}\int_{y}^{1}\left(x-y\right) dF_0\left(x\right)dG_0$ allows us to conclude the result. \end{proof}

\end{document}